\documentclass[]{interact}

\usepackage[pagewise]{lineno}
\usepackage{subfigure}

\usepackage{amsmath}
\usepackage{amsfonts}
\usepackage{amssymb}
\usepackage{graphicx}
\usepackage{mathrsfs}
\usepackage{color}

\usepackage{epstopdf}
\usepackage[caption=false]{subfig}

\usepackage[numbers,sort&compress]{natbib}
\bibpunct[, ]{[}{]}{,}{n}{,}{,}
\makeatletter
\def\NAT@def@citea{\def\@citea{\NAT@separator}}
\makeatother

\theoremstyle{plain}
\newtheorem{theorem}{Theorem}[section]
\newtheorem{lemma}[theorem]{Lemma}
\newtheorem{corollary}[theorem]{Corollary}
\newtheorem{proposition}[theorem]{Proposition}

\theoremstyle{definition}
\newtheorem{definition}[theorem]{Definition}

\theoremstyle{remark}
\newtheorem{remark}{Remark}

\newcommand{\ud}{\mathrm{d}}
\newcommand{\R}{\mathbb{R}}

\begin{document}

\title{On solutions of a partial integro-differential equation in Bessel potential spaces with applications in option pricing models}

\author{
\name{Jos\'{e} M. T. S. Cruz\textsuperscript{a}
and
Daniel \v{S}ev\v{c}ovi\v{c}\textsuperscript{b}\thanks{CONTACT D.~\v{S}ev\v{c}ovi\v{c}. Email: sevcovic@fmph.uniba.sk}
}
\affil{\textsuperscript{a}
ISEG, University of Lisbon, Rua de Quelhas 6, 1200-781 Lisbon, Portugal; 
\textsuperscript{b}
Comenius University in Bratislava, Mlynsk\'a dolina, 84248 Bratislava, Slovakia}
}

\maketitle

\begin{abstract}

In this paper we focus on  qualitative properties of solutions to a nonlocal nonlinear partial integro-differential equation (PIDE). Using the theory of abstract semilinear parabolic equations we prove existence and uniqueness of a solution in the scale of Bessel potential spaces. Our aim is to generalize known existence results for a wide class of L\'evy measures including with a strong singular kernel.

As an application we consider a class of PIDEs arising in the financial mathematics. The classical linear Black-Scholes model relies on several restrictive assumptions such as liquidity and completeness of the market. Relaxing the complete market hypothesis and assuming a L\'evy stochastic process dynamics for the underlying stock price process we obtain a model for pricing options by means of a PIDE. We investigate a model for pricing call and put options on underlying assets following a L\'evy stochastic process with jumps. We prove existence and uniqueness of solutions to the penalized PIDE representing approximation of the linear complementarity problem arising in pricing American style of options under L\'evy stochastic processes. We also present numerical results and comparison of option prices for various L\'evy stochastic processes modelling underlying asset dynamics. 

\end{abstract}


\section{Introduction}

In this paper, we analyze solutions to the semilinear parabolic partial integro-differential equation (PIDE) of the form: 
\begin{eqnarray}
\frac{\partial u}{\partial \tau}(\tau,x) &=& \frac{\sigma^2}{2}\frac{\partial^2 u}{\partial x^2}(\tau,x)  + \omega \frac{\partial u}{\partial x}(\tau,x) + g(\tau, u(\tau,x))
\nonumber
\\
&&  + \int_{\mathbb{R}}\left[ u(\tau, x+z)-u(\tau, x)- (e^z-1) \frac{\partial u}{\partial x}(\tau, x)\, \right] \nu(\ud z),
\label{PDE-u}
\\
&&u(0,x)=u_0(x),
\nonumber
\end{eqnarray}
$x\in \mathbb{R}, \tau\in(0,T)$, where $g$ is H\"older continuous in the $\tau$ variable and it is Lipschitz continuous in the $u$ variable. Here $\nu $  is a positive Radon measure on $\mathbb{R}$ such that $\int_{\mathbb{R}} \min(z^2,1) \nu (\ud z)<\infty$. 

Our purpose is to prove  existence and uniqueness of a solution to (\ref{PDE-u}) in the framework of Bessel potential spaces. These functional spaces represent a nested scale $\{X^\gamma\}_{\gamma\ge0}$ of Banach spaces such that
\[
X^1\equiv D(A) \hookrightarrow X^{\gamma_1} \hookrightarrow X^{\gamma_2} \hookrightarrow X^0\equiv X,
\]
for any $0\le \gamma_2\le \gamma_1\le 1$ where $A$ is a sectorial operator in the Banach space $X$ with a dense domain $D(A)\subset X$. For example, if $A=-\Delta$ is the Laplacian operator in $\R^n$ with the domain $D(A)\equiv W^{2,p}(\R^n)\subset X\equiv L^p(\R^n)$ then $X^\gamma$ is embedded in the Sobolev-Slobodecki space $W^{2\gamma,p}(\R^n)$ consisting of all functions  having $2\gamma$-fractional derivative belonging to the Lebesgue space $L^p(\R^n)$ of $p$-integrable functions (cf. \cite{Henry1981}). In this paper, our goal is to prove existence and uniqueness of solutions to (\ref{PDE-u}) for a general class of the so-called admissible activity L\'evy measures $\nu$ satisfying suitable growth conditions at $\pm\infty$ and the origin. 

A motivation for studying solutions of the PIDE  (\ref{PDE-u}) arises from financial modeling. In the last four decades, the Black-Scholes model and its various generalizations become popular in the financial industry because of their simplicity and possibility to price options by means of explicit analytic formulas. However, practical application of the classical linear Black-Scholes equation has serious drawbacks, e.g. evidence from the stock market indicating that this model is less realistic as it assumes that the market is liquid, complete and without transaction costs. Moreover, sample paths of a Brownian motion are continuous, but stock prices of a typical company usually suffer from sudden jumps on an intra-day scale, making the price trajectories discontinuous. In the classical Black-Scholes model, the logarithm of the price process has a normal distribution. However, the empirical distribution of stock returns exhibits fat tails. Furthermore, if we calibrate theoretical prices to the market prices, we realize that the implied volatility is neither constant as a function of strike nor as a function of time to maturity, contradicting thus assumptions of the Black-Scholes model. Several alternatives have been proposed in the literature for generalization of this model. The models with jumps can, at least in part, solve problems inherent to the Black-Scholes model. Jump--diffusion models also have an important role in derivative markets. In the classical Black-Scholes model the market is assumed to be complete, implying that every pay-off can be perfectly replicated. On the other hand, in jump--diffusion models there is no perfect hedge and this way options are not redundant. 

taking into account jumps in the underlying asset process, the price $V(t,S)$ of an option on the underlying asset with a price $S$ and time $t\in[0,T]$ is a solution to the following nonlocal nonlinear partial integro-differential equation: 
\begin{eqnarray}
\frac{\partial V}{\partial t} &+&\frac{\sigma^2}{2} S^2 \frac{\partial^2 V}{\partial S^2} + r S \frac{\partial V}{\partial S}-rV \nonumber
\\
&+& \int_{\mathbb{R}}\left[ V(t,Se^z)-V(t,S)-(e^z-1) S \frac{\partial V}{\partial S}(t,S) \right] \nu(\ud z)=0,
\label{V-eq}
\end{eqnarray} 
$S>0, t\in (0,T)$. Here $\sigma>0$ is the volatility of the underlying asset process $\{S_t\}_{t\ge0}$, $r\ge 0$ is the risk-less interest rate and $\nu$ is a L\'evy measure. A solution $V$ is subject to the terminal condition $V(S,T)=\Phi(S)$ where $\Phi$ represents the pay-off diagram of a plain vanilla option, i.e. $\Phi(S)= (S-K)^+$ for a call option, or $\Phi(S)= (K-S)^+$ for a put option, $K>0$ is the strike price. 

In the case when the L\'evy measure $\nu$ is defined through the Dirac function, i.e. $ \nu(\ud z) =\delta(z) \ud z$ or $\nu\equiv0$  the aforementioned nonlocal PIDE reduces to the classical linear Black-Scholes linear PDE:
\[
\frac{\partial V}{\partial t} +\frac{\sigma^2}{2} S^2 \frac{\partial^2 V}{\partial S^2} + r S \frac{\partial V}{\partial S}-rV =0.
\]

In the past years, existence results of PIDE (\ref{PDE-u}) have been intensively studied in the literature. In \cite{BL82} A. Bensoussan and J.-L. Lions  (see Theorem 3.3 and Theorem 8.1) and also M. G. Garroni and J. L. Menaldi (see \cite{GarrMenal2002}) investigated the existence and uniqueness of classical solutions for the case $\sigma>0$. In \cite{RMikuPraga93} Mikulevicius and Pragarauskas extended these results for the case $\sigma=0$. Furthermore, in \cite{RMikuPraga2014},\cite{RMikuPraga14} they investigated existence and uniqueness of classical solutions in H\"{o}lder and Sobolev spaces of the Cauchy problem to the partial-integro-differential equation of the order of kernel singularity up to the second order. Qualitative results using the notion of viscosity solutions were provided by M. Crandall and P.-L. Lions in \cite{CL92}. They were  generalized to PIDEs by Awatif  \cite{Sayahf2007} and Soner \cite{Soner1986} for the first order operators and by Alvarez and Tourin \cite{OlivierAgnes1996}, Barles {\it et al.}  \cite{GuyBarles1997}, and Pham \cite{Huyen1998} for the second order operators. In \cite{Maria11},\cite{Maria12} Mariani and SenGupta  proved existence of weak solutions of a generalized integro-differential equation using the Schaefer fixed point theorem. On other hand, in \cite{Amster12}, Amster {\it et al.} proved the existence of solutions using the method of upper and lower solutions in a general domain in the case of several assets and for the regime-switching jump-diffusion model in \cite{Ionut2012}. In \cite{NBS15},\cite{NBS17} Arregui et al. applied the theory of abstract parabolic equations in Banach spaces (cf. \cite{Henry1981}) for the proof of existence and uniqueness of solutions of a system of nonlinear PDEs for pricing of XVA derivatives. In the recent paper Cruz and \v{S}ev\v{c}ovi\v{c} \cite{NBS19} investigated a nonlinear extension of the option pricing PIDE model (\ref{V-eq}) from numerical point of view.

As a motivation we consider a model for pricing vanilla call and put options on underlying assets following L\'evy stochastic processes. Using the theory of abstract semilinear parabolic equations we prove existence and uniqueness of solutions in the Bessel potential space representing a fractional power space of the space of Lebesgue $p$-integrable functions with respect to the second order Laplace  differential operator. We generalize known existence results for a wider class of L\'evy measures having strong singular kernel with the third order of singularity. We also prove existence and uniqueness of solutions to the penalized PIDE representing approximation of the linear complementarity problem for a PIDE arising in pricing American style of options. 

The paper is organized as follows. In Section 2 we recall typical examples of L\'evy measures arising in the financial modelling of stochastic processes with random jumps. We introduce a notion of an admissible activity L\'evy measure. We show that this class of L\'evy measures includes jump-diffusion finite activity measures present in e.g. Merton's or Kou's double exponential models as well as infinite activity L\'evy measure appearing in e.g. Variance Gamma, Normal Inverse Gaussian or the so-called CGMY models. Section 3 is devoted to the proof of the main result on existence and uniqueness of solution to the PIDE (\ref{PDE-u}) in the framework of the Bessel potential spaces $X^\gamma$ representing the fractional power spaces of the Lebesgue space $L^{p}(\R)$ with respect to the second order Laplacian operator.  We follow the methodology of abstract semilinear parabolic equations developed by Henry in \cite{Henry1981}. First, we provide sufficient conditions guaranteeing existence and uniqueness of a solution to the PIDE (\ref{PDE-u}) in Bessel potential spaces. In Section 4 we investigate qualitative properties of solutions to a PIDE of the Black-Scholes type arising in pricing derivative securities on underlying assets following L\'evy processes. Section 5 is focused on application of the results for the nonlinear extension of the Black-Scholes PIDE for pricing American style of put options by the penalization method. 
Finally, in Section 6 we present results of a numerical solution to PIDE Variance Gamma and Merton's models.

\section{Preliminaries, definitions and motivation }

A stochastic process $\{X_{t}, t\ge 0\}$ is called a L\'{e}vy stochastic process if its characteristic function has the following L\'{e}vy-Khintchine representation 
$\mathbb{E}\left[e^{i y X_{t}}\right]=e^{t\phi(y)}$ with 
\begin{equation*}
\phi(y)=-\frac{\sigma^{2}}{2}y^{2}+i\omega y+\int_{-\infty}^{+\infty} \left(e^{i y z}-1-i y z1_{|z| \leq 1}\right) \nu  (\ud z),
\end{equation*} 
where $\sigma \geq 0$, $\omega \in \mathbb{R}$, and $\nu $  is a positive Radon measure on $\mathbb{R} \setminus \left\{0\right\}$ satisfying: 
\begin{equation}\int_{\mathbb{R}} \min(z^2,1) \nu (\ud z)<\infty, 
\label{condinttwo}
\end{equation} 
(cf. \cite{Sato99},\cite{ConTan03},\cite{App04}). 

\begin{definition}\label{def-admissiblemeasure}
A L\'evy measure $\nu$ is called an admissible activity L\'evy measure if  there exists a nonnegative measurable function $h$ such that $\nu(\ud z) = h(z) \ud z$ such that 
\begin{equation}
0 \le  h(z)\le C_0 |z|^{-\alpha}\left(e^{D^{-}z}1_{z\geq 0}+e^{D^{+}z}1_{z< 0}\right)e^{-\mu z^{2}},
\label{growth_measure}
\end{equation}
for any $z\in\R$ and the shape parameters $\alpha\geq 0$, $D^{\pm}\in\mathbb{R}$ and $\mu\geq 0$. Here $C_0>0$ is a constant. 
\end{definition}

The condition  (\ref{condinttwo}) is satisfied for any measure $\nu$ belonging to the class of admissible activity L\'evy measures with shape parameters $0\le \alpha<3$, and either $\mu>0$ or $\mu=0$ and $D^-<0<D^+$.

\subsection{Examples of admissible L\'evy measures arising in the financial modelling}

The class of admissible activity L\'{e}vy measures includes various measures often used in financial modelling of underlying stock dynamics with jumps. For example, in the context of financial modelling the first jump-diffusion model was proposed by Merton in \cite{Merton76}. Its L\'{e}vy measure is given by:
\begin{equation}
\nu(\ud z)=\lambda \frac{1}{\delta\sqrt{2\pi}}e^{-\frac{(z-m)^2}{2\delta^2}}\ud z\,, 
\label{merton-density}
\end{equation}
where $m\in\R, \lambda, \delta>0,$ are given parameters. 

Another popular model is the so-called double exponential model introduced by Kou in \cite{Kou2002}. In this model, the L\'evy measure $\nu$ is given by
\begin{equation}
\nu (\ud z)=\lambda \left( \theta \lambda^{+} e^{- \lambda^{+} z}1_{z\ge 0}+ (1-\theta) \lambda^{-} e^{\lambda^{-} z}1_{z<0}\right)\ud z,
\label{double-density}
\end{equation}
where $\lambda>0$ is the intensity of jumps, $\theta$ is the probability of occurrence of positive jumps and the parameters $\lambda^\pm> 0$ correspond to the level of the decay of distribution of  positive and negative jumps. It implies that the distribution of jumps is asymmetric and the tails of the distribution of returns are semi-heavy. Both Merton's as well as Kou's measure $\nu$ belong to the class of the so-called finite activity L\'evy measures, i.e. $\nu(\mathbb{R})=\int_{\mathbb{R}}\nu(\ud z)<\infty$ having a finite variation $\int_{|z|\leq 1}|z|\nu(\ud z)<\infty$. 

As an example of infinite activity L\'evy processes we can consider the Variance Gamma (see \cite{DPE98}), Normal Inverse Gaussian (NIG) and CGMY processes (see \cite{BARNIE01}).  
The Variance Gamma process is a process with infinite activity, $\nu(\mathbb{R})=\int_{\mathbb{R}}\nu(\ud z)=\infty$ and finite variation, $\int_{|z|\leq 1}|z|\nu(\ud z)<\infty$ where 
\begin{equation}
\nu (\ud z)= C_0 |z|^{-1}e^{Az-B|z|}\ud z. 
\label{vargamma-density}
\end{equation} 
Here the parameters $A,B>0$ depend on the volatility and drift of the Brownian motion, $C_0>0$,  and the variance of a subordinator (the Gamma process) (see \cite{ConTan03}). The measure $\nu$ is an admissible activity L\'evy measure with shape parameters $\mu=0$, $D^{+}=A+B>0$, $D^{-}=A-B<0$, and $\alpha=1$. The NIG process is a process of infinite activity and infinite variation with the following L\'{e}vy measure:
\begin{equation}
\nu(\ud z)=C |z|^{-1} e^{A z}K_{1}\left(B |z|\right)\ud z,
\label{nig-density}
\end{equation}
where $A,B>0$ have the same meaning as in the Variance Gamma process. Here $K_1$ is the modified Bessel function of the second kind (see \cite{ConTan03}). Since $K_{1}(x)\sim \sqrt{\pi/2} x^{-1/2} e^{-x}$ as  $x\rightarrow \infty$ and $K_{1}(x)\sim x^{-1}$ as  $x\rightarrow 0$ (see \cite{Abramowitz}) the measure $\nu$ is an admissible activity L\'evy measure with the  shape parameters  $\mu=0$, $D^{+}=A+B>0$, $D^{-}=A-B<0$, $\alpha=2$.
The Variance Gamma and NIG processes are special cases of generalized hyperbolic models.

Finally, the so-called CGMY distribution process introduced by Carr {\it et al.} in \cite{CGMY98},\cite{PD99} has four parameters $C,G,M$ and $Y$ with the the L\'{e}vy measure given by: 
\begin{equation}
\nu( \ud z)= C_0 |z|^{-1-Y} \left( e^{G z} 1_{z<0}+e^{-M z} 1_{z>0}\right) \ud z,
\label{CGMY-density}
\end{equation} 
where $C,G,M>0$ and $Y<2$. The parameter $C$ measures the overall level of activity. The parameters $G$ and $M$ are the left and right tail decay parameters, respectively. When $G=M$ the distribution is symmetric. The process has infinite activity and finite variation when $Y\in (0,1)$ and infinite variation for $Y\in[1,2)$. The measure $\nu$ is an admissible activity L\'evy  measure with the shape parameters $\mu=0, \alpha=1+Y<3$, and $D^{+}=G>0, D^{-}=-M<0$.

\section{Existence and uniqueness results}

The goal of this section is to prove the main result of the paper regarding existence and uniqueness of a solution to the linear and nonlinear PIDE for a wide class of admissible activity L\'evy measures. We can rewrite the PIDE (\ref{PDE-u}) in the abstract form as follows:
\begin{eqnarray}
&&\frac{\partial u}{\partial \tau} + A u = 
\omega \frac{\partial u}{\partial x} + f[u] + g(\tau, u), \quad x\in \mathbb{R}, \tau\in(0,T),
\label{problem_transformed}
\\
&&u(0,x)=u_0(x), x\in \mathbb{R},
\nonumber
\end{eqnarray}
where the linear operators $A$ and $f$ are defined by:
\begin{eqnarray}
A u &=& - \frac{\sigma^2}{2}\frac{\partial^2 u}{\partial x^2},
\label{Au_def}
\\
f[u](x) &=& 
\int_{\mathbb{R}}\left[ u(x+z)-u(x)- (e^z-1) \frac{\partial u}{\partial x}(x)\, \right] \nu(\ud z),
\label{functional_f_def}
\end{eqnarray}
and $g$ is a H\"older continuous mapping in the $\tau$ variable and it is Lipschitz continuous in the $u$ variable.

As a motivation for studying PIDE (\ref{problem_transformed}) we consider a model for pricing vanilla call and put options on underlying assets following L\'evy stochastic processes. The classical linear Black-Scholes equation can be transformed into equation v(\ref{problem_transformed}) where $f\equiv 0, g\equiv 0$. A nontrivial integral part $f[u]$ arises from a generalization of the Black-Scholes model to the case when the underlying asset price follows a stochastic L\'evy process with jumps (see Section 4). In Section 5 we will investigate equation  (\ref{problem_transformed}) with a nontrivial integral part $f[u]$ and a nonlinearity $g$ corresponding to the penalization function. The resulting PIDE of the form (\ref{problem_transformed}) represents an approximation of a solution to the partial integro-differential variational inequality arising in pricing American style of options.

In order to prove existence, continuation and uniqueness of a solution to the problem (\ref{problem_transformed}) we follow the qualitative theory of semilinear abstract parabolic equations developed by Henry in \cite{Henry1981}. First, we recall the concept of an analytic semigroup of linear operators and a sectorial operator in a Banach space. 

\begin{definition}\cite{Henry1981}
A family of bounded linear operators $\left\{S(t), t\geq 0\right\}$ in a Banach space $X$ is called an analytic semigroup if it satisfies the following conditions:
\begin{itemize}
\item[i)] $S(0)=I, S(t)S(s)=S(s)S(t)=S(t+s)$, for all $t,s\geq 0$;
\item[ii)] $S(t) u\rightarrow u$ when $t\rightarrow 0^{+}$ for all $u\in X$;
\item[iii)] $t\rightarrow S(t) u$ is a real analytic function on $0< t< \infty$ for each $u\in X$.
\end{itemize} 
The associated infinitesimal generator $A$ is defined as follows: $A u = \lim_{t\rightarrow 0^{+}} \frac{1}{t} (S(t)u-u)$ and its domain $D(A)\subseteq X$ consists of those elements $u\in X$ for which the limit exists in the space $X$. 
\end{definition}

\begin{definition}\cite{Henry1981}
Let $S_{a,\phi}=\left\{\lambda \in \mathbb{C}: \phi\leq \arg(\lambda-a)\leq 2\pi-\phi \right\}$ be a sector of complex numbers. A closed densely defined linear operator $A: D(A)\subset X \rightarrow X$ is called a sectorial operator if there exists a constant $M\geq 0$ such that $\Vert(A-\lambda )^{-1}\Vert\leq M/|\lambda -a|$ for all $\lambda\in S_{a,\phi} \subset \mathbb{C}\setminus\sigma(A)$.
\end{definition}

In what follows, we shall investigate the partial-integral differential equation (\ref{problem_transformed}) in the framework of the so-called Bessel potential spaces. These spaces represent natural extension of the classical Sobolev spaces $W^{k,p}(\mathbb{R})$ where the order $k$ may attain the discrete values only, i.e. the distributional derivatives up to the order $k$ belong to the Lebesgue space $L^p(\mathbb{R})$. Bessel potential spaces represent a continuous scale of fractional powers, and allow for a finer formulation of results in comparison to the classical Sobolev spaces $W^{k,p}(\mathbb{R}), k\in\mathbb{N}$.

It is well known that that if $A$ is a sectorial operator then $-A$ is an infinitesimal generator of an analytic semigroup $S(t)=\left\{e^{-A t}, t\geq 0\right\}$ (cf. \cite{Henry1981}). If $X$ is a Banach space then we can define a scale of fractional power spaces $\{X^{\gamma}\}_{\gamma\ge 0}$ in the following way:
\[
X^{\gamma}= D(A^{\gamma})= Range(A^{-\gamma})=\left\{u\in X:\  \exists \varphi\in X,  u=A^{-\gamma}\varphi\right\},
\]
where, for any $\gamma >0$, the operator $A^{-\gamma}$ is defined by virtue of the Gamma function, i.e. 
$A^{-\gamma}=\frac{1}{\Gamma(\gamma)}\int_{0}^{\infty} \xi^{\gamma-1}e^{-A \xi} \ud \xi$. The norm is defined as $\Vert u\Vert_{X^\gamma}=\Vert A^{\gamma}u\Vert_X=\Vert \varphi\Vert_X$. Note that $X^0=X$, $X^1=D(A)$, and $X^1\equiv D(A) \hookrightarrow X^{\gamma_1} \hookrightarrow X^{\gamma_2} \hookrightarrow X^0\equiv X$, for any $0\le \gamma_2\le \gamma_1\le 1$.

In what follows, by $G*\varphi$ we shall denote the convolution operator defined by $(G*\varphi)(x)=\int_{\mathbb{R}^n} G(x-y)\varphi(y)\ud{y}$. 

\begin{lemma}\cite[Section 1.6]{Henry1981}, \cite[Chapter 5]{Stein1970}
\label{A-sectorial}
The Laplace operator $-\Delta $ is sectorial in the Banach space $X=L^{p}(\mathbb{R}^{n})$ of Lebesgue $p$-integrable functions for any $p\ge 1$ and $n\ge 1$. Its domain $D(A)$ is embedded into the Sobolev space $W^{2,p}(\mathbb{R}^{n})$. The fractional power space $X^\gamma, \gamma>0,$ is the space of Bessel potentials:  $X^\gamma={\mathscr L}^p_{2\gamma}(\mathbb{R}^n):=\{ G_{2\gamma}*\varphi, \ \varphi\in L^p(\R^n)\}$ where 
\[
G_{2\gamma}(x) = \frac{(4\pi)^{-n/2}}{\Gamma(\gamma)}\int_0^\infty \xi^{-1+(2\gamma-n)/2} e^{-(\xi +\Vert x\Vert^2/(4\xi))}\ud \xi
\]
is the Bessel potential function. The norm of $u=G_{2\gamma}*\varphi$ is given by $\Vert u\Vert_{X^\gamma}=\Vert \varphi\Vert_{L^p}$. The fractional power space $X^\gamma$ is continuously embedded into the fractional  Sobolev-Slobodeckii space $W^{2\gamma,p}(\R^n)$.
\end{lemma}

\begin{remark}
Lemma~\ref{A-sectorial} was proven in \cite[Section 1.6]{Henry1981}, \cite[Chapter 5]{Stein1970}. The idea of the proof of sectoriality of the Laplace operator $-\Delta $  in the Banach space $X=L^p(\mathbb{R}^n)$ is based on estimation of the resolvent operator $(\lambda -\Delta)^{-1}$ in the $L^p$ norm. The rest of the proof of  Lemma~\ref{A-sectorial} is based on the analysis of the Fourier transform of the equation $(\lambda-\Delta)u=f$. The Fourier transform $\hat u$ of its solution is given by $\hat u(\xi) = (\lambda +|\xi|^2)^{-1}\hat f(\xi)$. The function $G_\alpha$ is then constructed  by means of the inverse Fourier transform of $\hat G_\alpha(\xi)  = (1+|\xi|^2)^{-\alpha/2}$, and $G_\alpha$ is given as in Lemma~\ref{A-sectorial}. For further details we refer the reader to Section 1.6 of \cite{Henry1981} and \cite[Chapter 5]{Stein1970}.
\end{remark}

\begin{lemma}\label{lemma-f} Assume $\nu$ is an admissible activity L\'evy measure with shape parameters $\alpha, D^\pm$, and $\mu$ where $\alpha<3$ and either $\mu>0, D^\pm\in\mathbb{R}$, or $\mu=0, D^- + 1<0<D^+$. Suppose that  $\gamma\ge 1/2$ and $\gamma>(\alpha-1)/2$. Then, for the mapping $f$ defined by (\ref{functional_f_def}),
there exists a constant $C>0$ such that, for any $u$ satisfying $\partial_x u \in X^{\gamma-1/2}$, the following estimate holds:
\[
\Vert f[u]\Vert_{L^p} \le C \Vert \partial_x u\Vert_{X^{\gamma-1/2}}.
\]
In particular, if $u\in X^\gamma$ we have $\Vert f[u]\Vert_{L^p} \le C \Vert u\Vert_{X^\gamma}$ and the mapping $f$ is a bounded linear operator from the fractional power space $X^\gamma$ into $X=L^p(\R)$. 
\end{lemma}

\noindent Proof. The mapping $f$ can be split as follows:  $f[u]=\tilde f[u] + \tilde\omega \partial_x u$ where
\begin{eqnarray*}
\tilde f[u](x) &=& \int_\R \left(u(x+z) - u(x) - z \frac{\partial u}{\partial x}(x) \right) \nu(\ud{z}),
\end{eqnarray*}
and $\tilde\omega=\int_\R \left(z - e^z+1 \right) \nu(\ud{z})$. Since $z - e^z+1=O(z^2)$ as $z\to0$, and
\[
0\le \nu(\ud{z}) = h(z) \ud{z} \le  |z|^{-\alpha} \tilde h(z) \ud{z},
\ \ 
\text{where}
\ \ 
\tilde h(z) = C_0 e^{-\mu z^2} \left(e^{D^-z} 1_{z\ge0} + e^{D^+z} 1_{z<0} \right),
\]
we have $\tilde\omega\in\R$ provided that $0\le \alpha<3$, and, either $\mu>0, D^\pm\in\R$, or $\mu=0$ and $D^- +1<0<D^+$.

First, we consider the case when $\gamma>1/2$. We shall prove boundedness of the second linear operator $\tilde f$. If $u$ is such that  $\partial_xu\in X^{\gamma-1/2}$ then there exists $\varphi\in X=L^p(\R)$ such that $\partial_x u = A^{-(2\gamma-1)/2}\varphi  = G_{2\gamma-1} * \varphi$ and 
\[
\Vert \partial_x u\Vert_{X^{\gamma-1/2}} = \Vert \varphi\Vert_X= \Vert \varphi\Vert_{L^p}.
\]
Hence, for any $x, \theta$, and $z$ we have
\[
\frac{\partial u}{\partial x}(x+\theta z)-  \frac{\partial u}{\partial x}(x)
= \left(G_{2\gamma-1}(x+\theta z - \cdot) - G_{2\gamma-1}(x - \cdot) \right)* \varphi(\cdot).
\]
Recall the following inequality for the convolution operator:
\[
\Vert G*\varphi\Vert_{L^p}\le \Vert G\Vert_{L^q} \Vert \varphi\Vert_{L^r},
\]
where  $p,q,r\ge 1$ and $1/p + 1 = 1/q + 1/r$ (see \cite[Section 1.6]{Henry1981}). In the special case when $q=1$ we have $\Vert G*\varphi\Vert_{L^p}\le \Vert G\Vert_{L^1} \Vert \varphi\Vert_{L^p}$. According to \cite[Chapter 5.4, Proposition 7]{Stein1970} we know that the modulus of continuity of the Bessel kernel function $G_{2\gamma-1}$  satisfies the estimate:
\[
\Vert G_{2\gamma-1}(\cdot + h) - G_{2\gamma-1}(\cdot) \Vert_{L^1}
\le C_1 |h|^{2\gamma-1},
\]
for any $h$ where $C_1>0$ is a constant. Therefore, for any $\theta, z\in\R$ we have
\begin{eqnarray*}
&&\int_\R \left|\frac{\partial u}{\partial x}(x+\theta z) -  \frac{\partial u}{\partial x}(x)\right|^p\ud{x}
= \Vert\left( G_{2\gamma-1}(\cdot + \theta z) - G_{2\gamma-1}(\cdot)\right)*\varphi\Vert_{L^p}^p 
\\
&& \le \Vert G_{2\gamma-1}(\cdot + \theta z) - G_{2\gamma-1}(\cdot)\Vert_{L^1}^p\Vert\varphi\Vert_{L^p}^p 
\le C_1^p |\theta z|^{(2\gamma-1)p} \Vert \partial_x u\Vert_{X^{\gamma-1/2}}^p.
\end{eqnarray*}
The latter inequality formally holds true also for the case $\gamma=1/2$ because 
\[
\int_\R \left|\frac{\partial u}{\partial x}(x+\theta z) -  \frac{\partial u}{\partial x}(x)\right|^p\ud{x}
\le 2^p \Vert \partial_x u\Vert_{L^p}^p = 2^p \Vert \partial_x u\Vert_{X^0}^p . 
\]
The rest of the proof of boundedness of the mapping $f$ holds for $\gamma> 1/2$ as well as $\gamma=1/2$.
As $u(x+z) - u(x) - z \frac{\partial u}{\partial x}(x)
=z\int_0^1 \frac{\partial u}{\partial x}(x+\theta z)-  \frac{\partial u}{\partial x}(x)\ud{\theta}$, we obtain
\begin{eqnarray*}
&& \int_\R |u(x+z) - u(x) - z \frac{\partial u}{\partial x}(x)|^p\ud{x}
=|z|^p \int_\R \left| \int_0^1 \frac{\partial u}{\partial x}(x+\theta z)-  \frac{\partial u}{\partial x}(x)\ud{\theta} \right|^p \ud{x}
\\
&& \le   |z|^p \int_0^1 \int_\R \left|  \frac{\partial u}{\partial x}(x+\theta z)-  \frac{\partial u}{\partial x}(x)\right|^p \ud{x} \ud{\theta} 
\le C_1^p |z|^{2\gamma p} \Vert \partial_x u\Vert_{X^{\gamma-1/2}}^p.
\end{eqnarray*}
Now, as $0\le \nu(\ud{z}) = h(z)\ud{z} \le  |z|^{-\alpha} \tilde h(z) \ud{z}=
(|z|^{-\beta} \tilde h(z)^\frac12 )\cdot (|z|^{\beta-\alpha} \tilde h(z)^\frac12 ) \ud{z}$, using the H\"older inequality with exponents $p,q$ such that $1/p+1/q=1$ we obtain
\begin{eqnarray*}
\Vert \tilde f[u]\Vert_{L^p}^p 
&=& \int_\R\left| \int_\R u(x+z) - u(x) - z \frac{\partial u}{\partial x}(x)\nu(\ud{z})\right |^p\ud{x}
\\
&\le& \int_\R \left|\int_\R \left| u(x+z) - u(x) - z \frac{\partial u}{\partial x}(x)\right| h(z)\ud{z}  \right|^p\ud{x}
\\
&\le& \int_\R \int_\R \left| u(x+z) - u(x) - z \frac{\partial u}{\partial x}(x)\right|^p |z|^{-\beta p} \tilde h(z)^{p/2} \ud{z} 
\\
&& \quad\times \left(\int_\R |z|^{(\beta-\alpha) q} \tilde h(z)^{q/2} \ud{z}\right)^{p/q} \ud{x}
\\
&=& \int_\R \left(\int_\R \left| u(x+z) - u(x) - z \frac{\partial u}{\partial x}(x)\right|^p \ud{x}\right) |z|^{-\beta p} \tilde h(z)^{p/2} \ud{z} 
\\
&& \quad\times \left(\int_\R |z|^{(\beta-\alpha) q} \tilde h(z)^{q/2} \ud{z}\right)^{p/q} 
\\
&\le & 
C_1^p  \Vert \partial_x u\Vert_{X^{\gamma-1/2}}^p
\int_\R |z|^{(2\gamma-\beta) p} \tilde h(z)^{p/2} \ud{z} \left(\int_\R |z|^{(\beta-\alpha) q} \tilde h(z)^{q/2} \ud{z}\right)^{p/q}.
\end{eqnarray*}
The integrals 
$C_2=\int_\R |z|^{(2\gamma-\beta) p} \tilde h(z)^{p/2} \ud{z}$ and 
$C_3=\int_\R |z|^{(\beta-\alpha) q} \tilde h(z)^{q/2} \ud{z}$
are finite provided that 
\[
(2\gamma-\beta) p > -1, \qquad (\beta-\alpha) q = (\beta-\alpha) \frac{p}{p-1} >-1,
\]
and $\mu>0, D^\pm\in\mathbb{R}$, or $\mu=0$ and $D^-<0<D^+$. 
The later inequalities are satisfied if there exists a parameter $\beta$ such that
\[
\alpha-1+1/p < \beta < 2\gamma+1/p.
\]
Such a choice of $\beta$ is possible because we assumed $\gamma>(\alpha-1)/2$. Hence there exists a constant $C>0$ such that $\Vert \tilde f[u]\Vert_{L^p} \le C \Vert \partial_x u\Vert_{X^{\gamma-1/2}}$ for any $u$ satisfying $\partial_x u\in X^{\gamma-1/2}$, as claimed. Due to the continuity of the embedding $X^{\gamma-1/2} \hookrightarrow X$ we have $\Vert f[u]\Vert_{L^p} = \Vert \tilde f[u] + \tilde\omega \partial_x u \Vert_{L^p} \le C \Vert \partial_x u\Vert_{X^{\gamma-1/2}} = C \Vert u\Vert_{X^\gamma}$ for any $u\in X^\gamma$ and $f$ is a bounded linear operator from $X^\gamma$ into $X=L^p$. \hfill$\diamondsuit$

Let us denote by $C([0,T],X^{\gamma})$ the Banach space of all continuous functions from the interval $[0,T]$ to $X^\gamma$ with the maximum norm $\Vert U(\cdot)\Vert_{C([0,T],X^{\gamma})}=\sup_{\tau\in[0,T]} \Vert U(\tau)\Vert_{X^\gamma}$. We recall the  well known result on existence and uniqueness of a solution to abstract parabolic equations in Banach spaces due to Henry \cite{Henry1981}.

\begin{proposition}\cite[Section 1]{Henry1981}
\label{semilinear_general_existence_result}
Suppose that a densely defined closed linear operator $-A$ is a generator of an analytic semigroup $\left\{e^{-At},t\geq 0\right\}$ in a Banach space $X$, $U_{0}\in X^{\gamma}$ where $0\leq \gamma <1$. Assume $F:[0,T]\times X^{\gamma}\to X$ and $h:(0,T]\to X$ are H\"older continuous mappings in the $\tau$ variable, $\int_0^T \Vert h(\tau)\Vert_X \ud x <\infty$, and $F$ is a Lipschitz continuous mapping in the $U$ variable. Then, there exists the unique solution $U\in C([0,T],X^{\gamma})$ of the following abstract semilinear evolution equation:
\begin{equation}
\frac{\partial U}{\partial \tau}+A U=F(\tau, U) +h(\tau), \qquad U(0)=U_{0}.
\label{semilinear_problem}
\end{equation}
Moreover, $\partial_\tau U(\tau) \in X, U(\tau)\in D(A)$ for any $\tau\in (0,T)$. 
\end{proposition}

\begin{remark}
By a solution to (\ref{semilinear_problem}) we mean a function $U\in C([0,T],X^{\gamma})$ satisfying  (\ref{semilinear_problem}) in the integral  (mild) sense, i.e. 
\[
U(\tau) = e^{-A \tau} U_0 + \int_0^\tau e^{-A (\tau-s)} (F(s, U(s)) + h(s) ) \ud{s} \ \hbox{for any}\ \tau\in[0,T].
\]
Recall that the key idea of the proof of Proposition~\ref{semilinear_general_existence_result} is based on the Banach fixed point argument combined with the decay estimate $\Vert e^{-A t}\Vert_{X^\gamma} = \Vert A^\gamma e^{-A t}\Vert_{X} \le M t^{-\gamma} e^{-at}$ of the norm of the semigroup $e^{-At}$ for any $t>0$.
\end{remark}

As a direct consequence of Proposition~\ref{semilinear_general_existence_result} and Lemma~\ref{lemma-f} we deduce the following result:
\begin{theorem}
\label{semilinear_existence_result}
Assume $\nu$ is an admissible activity L\'evy measure with the shape parameters $\alpha, D^\pm$ and $\mu$ where $\alpha<3$ and either $\mu>0, D^\pm\in\mathbb{R}$, or $\mu=0, D^- + 1<0<D^+$. Assume $\gamma\ge 1/2$ and $\gamma>(\alpha-1)/2$. Suppose that the function $g(\tau,u)$ is H\"older continuous in the $\tau$ variable and Lipschitz continuous in the $u$ variable. Then for any $u_0\in X^\gamma$ and $T>0$ there exists the unique solution $u\in C([0,T],X^{\gamma})$ to PIDE (\ref{PDE-u}).
\end{theorem}

\section{The Black-Scholes PIDE model}
In this section, our purpose is to investigate properties of solutions to a PIDE generalizing the Black-Scholes model. An important definition concerning this generalization is definition of a L\'{e}vy measure of a given process $X_t$.
The measure $\nu(A)$ of a Borel set $A\subseteq \R$ is defined by: 
\begin{equation}
\nu\left(A\right)=\mathbb{E}\left[\# \left\{t \in \left[0,1\right]:\Delta X_{t} \in A \right\}\right]=
\frac{1}{T}\mathbb{E}\left[\# \left\{t \in \left[0,T\right]:\Delta X_{t} \in A \right\}\right].
\label{eq:measuredef}
\end{equation} 
It gives the mean number, per unit of time, of jumps of $X_t, t\ge0,$ whose amplitude belongs to the set $A$ (see \cite{ConTan03}). 

For the underlying asset price dynamics  we will suppose that $S_t, t\ge 0,$follows the geometric L\'evy proces, i.e. $S_t=e^{X_t}$ where $X_t,t\ge0,$ is a L\'evy process. Then it is well known (cf. \cite{ConTan03},\cite{NBS19}) that the price of a contingent claim in the presence of jumps is given by a solution $V(t,S)$ of the following partial integro-differential equation:
\begin{eqnarray}
\frac{\partial V}{\partial t} &+&\frac{\sigma^2}{2} S^2 \frac{\partial^2 V}{\partial S^2} + r S \frac{\partial V}{\partial S}-rV \nonumber
\\
&+&\int_{\mathbb{R}}\left[ V(t,Se^z)-V(t,S)-(e^z-1) S \frac{\partial V}{\partial S}(t,S) \right] \nu(\ud z)=0,
\label{PDE-S}
\\
&&V(T,S)=\Phi(S), \quad S>0, t\in[0,T).
\nonumber
\end{eqnarray} 
Here $\Phi$ is the pay-off diagram of a plain vanilla option. For example,  $\Phi(S)=(S-K)^+$ for a call option, or $\Phi(S)=(K-S)^+$ for a put option where $K>0$ is the strike price. Here and after we shall denote by $a^+=\max(a,0)$ and $a^-=\min(a,0)$ the positive and negative parts of a real number $a$, respectively.

If we consider the following change of variables $V(t,S)=e^{-r\tau}u(\tau,x)$ where $\tau=T-t$, $x=\ln(\frac{S}{K})$ then we obtain the following PIDE for the function  $u(\tau,x)$: 

\begin{eqnarray}
\frac{\partial u}{\partial \tau}
&=&\frac{\sigma^2}{2} \frac{\partial^2 u}{\partial x^2} 
+ \left(r-\frac{1}{2}\sigma^2\right)\frac{\partial u}{\partial x}
\label{PDE-u-BS}
\\
&& + \int_{\mathbb{R}} \left[
u(\tau,x+z)-u(\tau,x)-\left(e^{z}-1\right)\frac{\partial u}{\partial x}(\tau,x)
\right]\nu(\ud z), 
\nonumber
\\
u(0,x)&=&\Phi(K e^{x}), \quad  x\in \mathbb{R}, \tau\in(0,T).
\nonumber
\end{eqnarray}

Unfortunately, the initial  condition $u(0,x)= \Phi(Ke^x)$ does not belong to the Banach space $X$ for both call and put option pay-off diagrams $\Phi$, i.e. $\Phi(S)=(S-K)^+$ and $\Phi(S)=(K-S)^+$. The idea how to formulate existence and uniqueness of a solution to the PIDE  (\ref{PDE-u-BS}) is based on the idea of shifting the solution $u$ by $u_{BS}$ where the function $u_{BS}(\tau,x) =e^{r\tau} V_{BS}(T-\tau, Ke^x)$ corresponds to transformation of the classical solution $V_{BS}$ to the linear Black-Scholes equation without PIDE part, i.e.
\begin{eqnarray*}
&&\frac{\partial V_{BS}}{\partial t} +\frac{\sigma^2}{2} S^2 \frac{\partial^2 V_{BS}}{\partial S^2} + r S \frac{\partial V_{BS}}{\partial S}-rV_{BS} =0, 
\\
&&V_{BS}(T,S)=\Phi(S).
\nonumber
\end{eqnarray*} 
Recall that the solution $V_{BS}$ for a call or put option can be expressed explicitly: 
\begin{eqnarray*}
V^{call}_{BS}(t,S) &=& S N(d_1) - K e^{-r(T-t)} N(d_2), 
\\
V^{put}_{BS}(t,S)  &=& K e^{-r(T-t)} N(-d_2) - S N(-d_1), 
\end{eqnarray*}
where
\[
d_{1,2} = \frac{\ln(S/K) + (r\pm\sigma^2/2)(T-t)}{\sigma\sqrt{T-t}},
\ \ \hbox{and} \ \  N(d)=\int_{-\infty}^d \frac{e^{-\xi^2/2}}{\sqrt{2\pi}} \ud \xi
\]
is the cumulative distribution function of the normal distribution (cf. \cite{NBS5}). Furthermore, the transformed function $u_{BS}$ is a solution to the linear parabolic PDE:
\begin{eqnarray}
&& \frac{\partial u_{BS}}{\partial \tau}
=\frac{\sigma^2}{2} \frac{\partial^2 u_{BS}}{\partial x^2} 
+ \left(r-\frac{1}{2}\sigma^2\right)\frac{\partial u_{BS}}{\partial x},
\label{PDE-uBS}
\\
&&u_{BS}(0,x)=\Phi(K e^{x}), \quad \tau\in(0,T), x\in \mathbb{R},
\nonumber
\end{eqnarray}
where $\Phi(Ke^x)=K(e^x-1)^+$ for the  call option and  $\Phi(Ke^x)=K(1-e^x)^+$ for the put option.

In what follows, we shall provide important estimates for the function $f[u_{BS}]$.

\begin{lemma}\label{lemmafuBS}
Suppose that $\nu$ is  an admissible activity L\'evy measure $\nu$ with the shape parameters $\alpha, D^\pm$, and $\mu$ where $\alpha<3$ and either $\mu>0, D^\pm\in\mathbb{R}$, or $\mu=0, D^- +1<0<D^+$. Suppose that 
$\frac12 \le \gamma<1$ and $\frac{\alpha-1}{2}<\gamma < \frac{p+1}{2p}\le1$. Then there exists a constant $C_0>0$ depending on the parameters $p,\sigma,r,T,K$ only, and such that the function $f[u_{BS}(\tau, \cdot)]$ satisfies the following estimates:
\begin{eqnarray*}
&&\Vert f[u_{BS}(\tau, \cdot)] \Vert_{L^p} 
\le C_0 \tau^{-(2\gamma-1)\left(\frac{1}{2} - \frac{1}{2p}\right)},\qquad 0<\tau\le T,
\\
&&\Vert f[\partial_\tau u_{BS}(\tau, \cdot)] \Vert_{L^p} 
\le C_0 \tau^{-\gamma-\frac{1}{2} + \frac{1}{2p}}, \qquad 0<\tau\le T,
\\
&&\Vert f[u_{BS}(\tau_1, \cdot)] - f[u_{BS}(\tau_2, \cdot)] \Vert_{L^p} 
\le C_0 |\tau_1-\tau_2|^{-\gamma +\frac{p+1}{2p} }, \qquad 0<\tau_1,\tau_2\le T.
\end{eqnarray*}

\end{lemma}

\noindent Proof. First, we consider the case of a call option, i.e. $u_{BS}=u^{call}_{BS}$ with $u_{BS}(0,x)=\Phi(K e^x) = K (e^x - 1)^+$.  It is important to emphasize that $f[e^x]=0$. Hence
\[
f[u_{BS}] = f[u_{BS} - K e^{r\tau+x}], \quad \hbox{and}\ \ 
\partial_\tau f[u_{BS}] = f[\partial_\tau(u_{BS} - K e^{r\tau+x})].
\]
In what follows, we shall denote by $C_0$ any generic positive constant depending on the parameters $p,\sigma,r,T,K$ only.  
With regard to Lemma~\ref{lemma-f} we shall estimate the $X^{\gamma-1/2}$ norm of the function $v$:
\begin{equation}
v(\tau,x)= \partial_x \left(u_{BS}(\tau,x) - K e^{r\tau+x}\right)
= K e^{r\tau+x}( N(d_1(\tau,x)) -1),
\label{vtaux}
\end{equation}
where $d_1(\tau,x) = \left(x+(r+\sigma^2/2)\tau\right)/(\sigma\sqrt{\tau})$. 
In the case of a put option we have
\[
\partial_x u^{put}_{BS}(\tau,x) =  -K e^{r\tau +x} N(-d_1(\tau,x))
= -K e^{r\tau +x} (1-N(d_1(\tau,x))) = v(\tau, x).
\]
Hence the proof of the statement of lemma for the case of a put option is essentially the same as the following argument for a call option. 

Using integration by parts and substitution $\xi=d_1(\tau,x)$, we obtain  
\begin{eqnarray*}
\Vert v(\tau,\cdot)\Vert_{L^p}^p 
&=& K^p e^{p r \tau} \int_{-\infty}^\infty e^{p x} (1-N(d_1))^p \ud x
\\
&\le& 
K^p e^{p r \tau} \int_{-\infty}^\infty e^{p x} (1-N(d_1)) \ud x
= K^p e^{p r \tau} \int_{-\infty}^\infty \frac{e^{p x}}{p} \frac{e^{-d_1^2/2}}{\sqrt{2\pi}} \frac{1}{\sigma\sqrt{\tau}} \ud x 
\\
&=& 
K^p e^{p r \tau} \int_{-\infty}^\infty \frac{e^{p \sigma\sqrt{\tau}\xi - p(r+\sigma^2/2)\tau}}{p} \frac{e^{-\xi^2/2}}{\sqrt{2\pi}} \ud \xi 
= \frac{1}{p} K^p e^{p(p-1)\tau \sigma^2/2}.
\end{eqnarray*}
Thus $\Vert v(\tau,\cdot)\Vert_{L^p} \le  p^{-1/p} K e^{(p-1) T \sigma^2/2}\equiv C_0$ for any $0<\tau\le T$.

As $\partial_x v = v + w$ where 
\[
w = K e^{ r \tau+x} N'(d_1) \frac{1}{\sigma\sqrt{\tau}}
= K e^{ r \tau+x} \frac{e^{-d_1^2/2}}{\sigma\sqrt{2\pi\tau}}.
\]
we obtain 
\begin{eqnarray}
\Vert w(\tau,\cdot)d_1(\tau,\cdot)^k\Vert_{L^p}^p 
&=& \frac{K^p e^{p r \tau}}{(\sigma\sqrt{2\pi\tau})^{p-1}} \int_{-\infty}^\infty e^{p x} \frac{e^{-p d_1^2/2} |d_1|^{p k}} {\sigma\sqrt{2\pi\tau}} \ud x \nonumber
\\
&=& 
\frac{K^p e^{p r \tau}}{(\sigma\sqrt{2\pi\tau})^{p-1}} \int_{-\infty}^\infty
e^{p \sigma\sqrt{\tau}\xi - p(r+\sigma^2/2)\tau} \frac{e^{-\xi^2/2} |\xi|^{p k}}{\sqrt{2\pi}} \ud \xi 
\label{wdineq}
\\
&\le& C_0^p \tau^{-\frac{p-1}{2}} \nonumber
\end{eqnarray}
for $k=0,1,2$. Applying (\ref{wdineq}) with $k=0$ we obtain $\Vert w(\tau,\cdot)\Vert_{L^p} \le C_0 \tau^{-\frac{1}{2}+\frac{1}{2p}}$. As a consequence, $\Vert v(\tau,\cdot)\Vert_{W^{1,p}} \le C_0 \tau^{-\frac{1}{2} +\frac{1}{2p}}$. Since the Bessel potential space ${\mathscr L}^p_{2\gamma-1}$ is an interpolation space between ${\mathscr L}^p_0=L^p$ and ${\mathscr L}^p_1=W^{1,p}$ using the Gagliardo-Nirenberg interpolation inequality
\[
\Vert v\Vert_{X^{\gamma-1/2} }
\equiv \Vert v\Vert_{{\mathscr L}^p_{2\gamma-1}}
\le C_0 \Vert v\Vert_{L^p}^\theta \Vert v\Vert_{W^{1,p}}^{1-\theta}, \quad \hbox{where}\ 
2\gamma-1 = 0\cdot\theta + 1\cdot (1-\theta), 
\]
(cf. \cite[Section 1.6]{Henry1981}) and applying Lemma~\ref{lemma-f} we obtain 
\[
\Vert f[u_{BS}(\tau, \cdot)] \Vert_{L^p} 
\le C \Vert v(\tau,\cdot) \Vert_{X^{\gamma-1/2}}
\le  C_0 \tau^{-(2\gamma-1)\left(\frac{1}{2} - \frac{1}{2p}\right)},\qquad 0<\tau\le T,
\]
as claimed. 

In order to prove the remaining estimates, let us estimate the norm $\Vert \partial_\tau v(\tau,\cdot)\Vert_{X^{\gamma-1/2}}$. As $\partial_\tau d_1 = - \tau^{-3/2} x/(2\sigma) + \tau^{-1/2}(r+\sigma^2/2)/(2\sigma)=- \tau^{-1} d_1/2 + \tau^{-1/2}(r+\sigma^2/2)/\sigma$ we have 
\[
\partial_\tau v = r v + K e^{r\tau +x} N'(d_1)\partial_\tau d_1
=
r v + w ( -\tau^{-1/2} \sigma d_1/2 +r+\sigma^2/2 ).
\]
Using estimate (\ref{wdineq}) with $k=0,1$ we obtain 
\[
\Vert \partial_\tau v(\tau,\cdot)\Vert_{L^p} 
\le C_0 \tau^{-1 +\frac{1}{2p}},\qquad 0<\tau\le T.
\]
To estimate the $W^{1,p}$ norm of $\partial_\tau v$ we recall that $\partial_x v = v + w$. Thus 
\begin{eqnarray*}
\partial_x \partial_\tau v &=& \partial_\tau v + \partial_\tau w
=\partial_\tau v + r w 
+ K e^{r\tau +x} \left( \frac{N''(d_1)}{\sigma\sqrt{\tau}}\partial_\tau d_1 -\frac{N'(d_1)}{2\sigma\tau^{3/2}} \right)
\\
&=&
\partial_\tau v + r w + w \left(-d_1 \partial_\tau d_1  - \tau^{-1}/2 \right)
\\
&=& \partial_\tau v + r w + w \left(d_1^2 \tau^{-1}/2  - \tau^{-1}/2
-\tau^{-1/2} d_1 (r+\sigma^2/2)/\sigma  \right),
\end{eqnarray*}
as $N''(d_1)=-d_1 N'(d_1)$. Using estimate (\ref{wdineq}) with $k=0,1,2$, we obtain 
\[
\Vert \partial_\tau v(\tau,\cdot)\Vert_{W^{1,p}} 
\le C_0 \tau^{-\frac{3}{2} +\frac{1}{2p}},\qquad 0<\tau\le T.
\]
Again, using the Gagliardo-Nirenberg interpolation inequality
\[
\Vert \partial_\tau v\Vert_{X^{\gamma-1/2} }
\equiv \Vert \partial_\tau v\Vert_{{\mathscr L}^p_{2\gamma-1}}
\le C_0 \Vert \partial_\tau v\Vert_{L^p}^\theta \Vert \partial_\tau v\Vert_{W^{1,p}}^{1-\theta}, \quad \hbox{where}\ 
2\gamma-1 = 0\cdot\theta + 1\cdot (1-\theta)
\]
and applying Lemma~\ref{lemma-f} we obtain 
\[
\Vert \partial_\tau f[u_{BS}(\tau, \cdot)] \Vert_{L^p} 
\le C \Vert \partial_\tau v(\tau,\cdot) \Vert_{X^{\gamma-1/2}}
\le  C_0 \tau^{-\gamma-\frac{1}{2} + \frac{1}{2p}},\qquad 0<\tau\le T,
\]
as claimed in the second statement of lemma.

Finally, 
\begin{eqnarray*}
&&\Vert f[u_{BS}(\tau_1, \cdot)] - f[u_{BS}(\tau_2, \cdot)] \Vert_{L^p}
=
\Vert \int_{\tau_1}^{\tau_2} \partial_\tau f[u_{BS}(\tau, \cdot)] d\tau  \Vert_{L^p}
\\
&&\le \left|\int_{\tau_1}^{\tau_2} \Vert  \partial_\tau f[u_{BS}(\tau, \cdot)]   \Vert_{L^p} d\tau \right|
\le C_0 |\tau_1-\tau_2|^{-\gamma +\frac{p+1}{2p}}, \qquad 0<\tau_1,\tau_2\le T,
\end{eqnarray*}
and the function $f[u_{BS}(\tau,\cdot)]$ is H\"older continuous with the H\"older exponent  $-\gamma +\frac{p+1}{2p}>0$. The proof of lemma follows. \hfill$\diamondsuit$

\medskip

Combining the previous Lemmas~\ref{lemma-f}, \ref{lemmafuBS}, sectoriality of the operator $A=-\partial^2_x$ in $X=L^p(\mathbb{R})$ (see Lemma~\ref{A-sectorial}), and Proposition~\ref{semilinear_existence_result} we obtain the following existence and uniqueness result for the linear PIDE (\ref{PDE-u-BS}), and, consequently, for the linear option pricing model (\ref{PDE-S}):

\begin{theorem}\label{existence_linear_PIDE}
Assume $\nu$ is an admissible activity L\'evy  measure with the shape parameters $\alpha<3$ and either $\mu>0, D^\pm\in \mathbb{R}$, or $\mu=0$ and $D^- +1<0<D^+$. Let $X^\gamma={\mathscr L}^p_{2\gamma}(\mathbb{R})$ be the space of Bessel potentials where  $\frac12 \le \gamma<1$ and $\frac{\alpha-1}{2}<\gamma < \frac{p+1}{2p}$.

Then, for any $T>0$,  the linear PIDE (\ref{PDE-u-BS}) has the unique solution $u$ such that the difference $U=u-u_{BS}$ satisfies $U\in C([0,T],X^{\gamma})$. Moreover, $U(\tau,\cdot)\in X^1 =  {\mathscr L}^p_{2}(\mathbb{R})\subseteq W^{2,p}(\mathbb{R})$ and $\partial_\tau U(\tau, \cdot)\in X=L^p(\R)$ for any $\tau\in (0,T)$. 
\end{theorem}

\noindent Proof. Since the Black-Scholes solution $u_{BS}$ solves the linear PDE (\ref{PDE-uBS}) the difference $U=u-u_{BS}$ of a solution $u$ to (\ref{PDE-u-BS}) and $u_{BS}$ satisfies the PIDE:
\begin{eqnarray*}
\frac{\partial U}{\partial \tau}
&=&\frac{\sigma^2}{2} \frac{\partial^2 U}{\partial x^2} 
+ \left(r-\frac{1}{2}\sigma^2\right)\frac{\partial U}{\partial x} + f[U] + f[u_{BS}], 
\nonumber
\\
U(0,x)&=& 0, \quad x\in \mathbb{R}, \tau\in(0,T).
\nonumber
\end{eqnarray*}
This PIDE equation can be rewritten in the abstract form:
\begin{equation}
\frac{\partial U}{\partial \tau} + A U = F(U) + h(\tau),
\quad U(0)=0,
\label{problem_transformed-shifted}
\end{equation}
where the linear operators $A$ and $f$ were defined in (\ref{Au_def}) and (\ref{functional_f_def}). The functions $F=F(U)$ and $h=h(\tau)$, $F:X^\gamma\to X$, $h:(0,T]\to X$ are defined as follows:
\[
F(U) = (r-\sigma^2/2)\frac{\partial U}{\partial x} + f[U],
\qquad h(\tau)=f[u_{BS}(\tau,\cdot)].
\]
With regard to Lemma~\ref{lemma-f},  $F$ is a bounded linear mapping, and, consequently Lipschitz continuous from the space $X^\gamma$ into $X$ provided that $\gamma\ge 1/2$ and $\gamma>(\alpha-1)/2$. 

Taking into account Lemma~\ref{lemmafuBS} we obtain 
\[
\Vert h(\tau_1) - h(\tau_2) \Vert_{L^p} 
= \Vert f[u_{BS}(\tau_1, \cdot)] - f[u_{BS}(\tau_2, \cdot)] \Vert_{L^p} 
\le C_0 |\tau_1-\tau_2|^{-\gamma +\frac{p+1}{2p}}, 
\]
for any $0<\tau_1,\tau_2\le T$. Since $\gamma<\frac{p+1}{2p}$ the mapping $h:[0,T]\to X\equiv L^p(\mathbb{R})$ is H\"older continuous. Moreover, 
\[
\int_0^T \Vert h(\tau) \Vert_{L^p} d\tau=
\int_0^T \Vert f[u_{BS}(\tau, \cdot)] \Vert_{L^p}d\tau 
\le C_0 \int_0^T\tau^{-(2\gamma-1)\left(\frac{1}{2} - \frac{1}{2p}\right)}d\tau <\infty,
\]
because $(2\gamma-1)\left(\frac{1}{2} - \frac{1}{2p}\right)<1$. The rest of the proof now follows from Theorem~\ref{semilinear_existence_result}. 
\hfill $\diamondsuit$

\bigskip

The following corollary is a consequence of embedding of the Bessel potential space into the space of H\"older continuous functions. 

\begin{corollary}
Suppose that an admissible activity L\'evy  measure $\nu$ fulfills assumptions of Theorem~\ref{existence_linear_PIDE}. 
Then, for any $T>0$, the linear PIDE (\ref{PDE-u-BS}) has the unique solution $u\in C([0,T], C^\kappa_{loc}(\mathbb{R}))$, with the H\"older exponent $\kappa>0$ satisfying  $\alpha-1-1/p < \kappa<1$. 
\end{corollary}

\noindent Proof. Recall continuity of the embedding 
\[
X^\gamma = {\mathscr L}^p_{2\gamma}(\mathbb{R}) \hookrightarrow C^\kappa_{loc}(\mathbb{R}),
\]
where $\kappa= 2\gamma -1/p$ (cf. \cite[Section 1.6]{Henry1981}), i.e.  $\gamma=\kappa/2 + 1/(2p)$. Now, there exists $1/2\le \gamma<1$ such that $\frac{\alpha-1}{2}<\gamma < \frac{p+1}{2p}$ if and only if $\alpha-1-1/p < \kappa<1$, as claimed. Therefore $U=u-u_{BS}$ belongs to $C([0,T], C^\kappa_{loc}(\mathbb{R}))$. 

The solution $u_{BS}=u_{BS}(\tau,x)$ is a real analytic function in the $\tau$ and $x$ variables for any $\tau>0$ and $x\in\mathbb{R}$. As $u_{BS}(0,x)$ represents the transformed call or put payoff diagram we have $u_{BS}=u_{BS}(0,x)$ is locally Lipschitz continuous in the $x$ variable. Hence $u_{BS}\in C([0,T], C^\kappa_{loc}(\mathbb{R}))$. Therefore the solution $u=U+u_{BS}$ to the linear PIDE (\ref{PDE-u-BS}) belongs to $C([0,T], C^\kappa_{loc}(\mathbb{R}))$, as claimed. \hfill $\diamondsuit$

\bigskip

\begin{remark}
Our method of the proof of existence and uniqueness of solutions to PIDEs can be extended to the multidimensional case in which the underlying fractional power space is $X^\gamma = {\mathscr L}^p_{2\gamma}(\mathbb{R}^n), n>1$. Recently, SenGupta, Wilson and Nganje \cite{SenGupta2019} studied a two factor Barndorff-Nielsen and Shephard model ($n=2$) with stochastic volatility in which both the underlying asset price $S$ and the variance $\sigma^2$ follow two finite activity admissible activity L\'evy procesess with a shape parameter $\alpha<3$. Their model can be applied for construction of an optimal hedging strategy for oil extraction that is benefiting from fracking technology. 
\end{remark}

\begin{remark}
The conditions  $\frac12 \le \gamma<1$ and $\frac{\alpha-1}{2}<\gamma < \frac{p+1}{2p}$ are fulfilled for a power $p\ge 1$ provided that either $\alpha\in [0,2]$ and $p\ge1$, or $\alpha\in (2,3)$ and $1\le p < 1/(\alpha-2)$. It means that if the L\'evy measure $\nu$ has a strong singularity of the order $\alpha\in(2,3)$ at the origin then we can find a solution in the framework of fractional power spaces of the Banach space $X=L^p(\R)$ where $p$ is limited by the order $\alpha$. The advantage of the choice of the Bessel potential space $X^\gamma = {\mathscr L}^p_{2\gamma}(\mathbb{R}), 1/2\le\gamma<1$, consists in the fact that we can prove existence and uniqueness of solutions in the phase space $X^\gamma$ for the case of stronger singularities with the order of singularity $\alpha$ up to $3$. The usual choice of the Sobolev space $X^{1/2}=W^{1,p}(\R)$ leads to the restriction of the order $\alpha$ of the singularity to $\alpha<2$. 
\end{remark}

\section{Existence results for nonlinear PIDE option pricing models}

In this section we present an application of the general existence and uniqueness result for the penalized version of the PIDE for solving the linear complementarity problem arising in pricing American style of a put option on an underlying asset following L\'evy stochastic process.

In \cite{BL82} Bensoussan and Lions characterized price of a put option in terms of a solution of a system of partial-integro differential inequalities (see also \cite{Lamberton2007}). In \cite{Wang2006} and \cite{Wang2017} Wang {\it et al.} investigated a penalty method for solving a linear complementarity problem using a power penalty term for the case without jumps in the underlying asset dynamics. In \cite{Donny2015} Lesman and Wang proposed a  power penalty method for solving the free boundary problem for pricing American options under transaction costs. Penalty methods for American option pricing under stochastic volatility models are studied in the paper \cite{Zvan1998} by Zvan, Forsyth and Vetzal. In \cite{Halluin2004} d'Halluin, Forsyth, and Labahn investigated a penalty method for American options on jump diffusion underlying processes.

Recall that American style options can be exercised anytime before the maturity time $T$. In the case of an American put option the state space $\{(t,S),\, t\in[0,T], S>0\}$ can be divided into the so-called early exercise region $\mathcal E$ and continuation region $\mathcal C$ where the put option should be exercised and hold, respectively. These regions are separated by the early exercise boundary defined by a function $t\mapsto S_f(t)$, such that $0<S_f(t)\le K$, and
\[
{\mathcal E} = \{ (t,S),\, t\in[0,T], 0<S\le S_f(t)\},
\quad {\mathcal C} = \{ (t,S),\, t\in[0,T], S_f(t)< S\}. 
\]
We refer the reader to papers \cite{NBS5}, \cite{SSC1999}, \cite{LS2010}, \cite{ZHU2006} for an overview of qualitative properties of the early exercise boundary for the case of pricing American style of put options for the Black-Scholes PDE with no integral part.

In the continuation region $\mathcal C$ the put option price is strictly greater than the pay-off diagram, i.e. $V(t,S)>\Phi(S)=(K-S)^+$ for $S_f(t)<S$. In the exercise region $\mathcal E$ the put option price is given by its pay-off diagram, i.e. $V(t,S) = \Phi(S)=(K-S)^+$. Moreover, the put option price $V(t,S)$ is a decreasing function in the $S$ variable. Hence in the exercise region where $0<S<S_f(t)\le K$, for the price $V(t,S)=K-S$ we obtain 
\begin{eqnarray*}
\frac{\partial V}{\partial t} + L^{S}[V]
&\equiv& 
\frac{\partial V}{\partial t} + 
\frac{\sigma^{2}}{2}S^{2}\frac{\partial^{2}V}{\partial S^{2}}
+ r S \frac{\partial V}{\partial S}-rV
\\
&& +\int_{\mathbb{R}}\left[V(t,Se^{y})-V(t,S)-S\left(e^{y}-1\right)
\frac{\partial V}{\partial S}(t,S) \right]\nu(\ud y)
\\
&=& 
-r K + 
\int_{-\infty}^0\left[V(t,Se^{y})-(K-S)- S\left(e^{y}-1\right)(-1) \right]\nu(\ud y)
\\
&& + \int_0^\infty\left[V(t,Se^{y})-(K-S)- S\left(e^{y}-1\right)(-1) \right]\nu(\ud y)
\\
&=& 
-r K  + \int_0^\infty\left[V(t,Se^{y})-(K-S) + S\left(e^{y}-1\right) \right]\nu(\ud y)
\\
&\le & 
-r K  + S \int_0^\infty\left(e^{y}-1\right)\nu(\ud y)
\end{eqnarray*}
because $S\mapsto V(t,S)$ is a decreasing function, and thus  $V(t,Se^{y})\le V(t,S)= K-S$ for $y\ge 0$, and $V(t,Se^{y})= K-Se^{y}$ for $y\le 0$. 

Let us assume that the admissible activity L\'evy measure $\nu$ satisfies the inequality:
\begin{equation}
    \int_0^\infty\left(e^{y}-1\right)\nu(\ud y) \le r.
\label{nuassumption}
\end{equation}
Then the price $V(t,S)$ of an American put option satisfies the inequality $\partial_t V(t,S) + L^{S}[V](t,S)\le 0$ for $0<S\le S_f(t)\le K$, i.e. for  $(t,S)\in{\mathcal E}$. On the other hand, for   $(t,S)\in{\mathcal C}$ the price $V(t,S)$ is obtained from the Black-Scholes PIDE equation $\partial_t V(t,S) + L^{S}[V](t,S)= 0$.

In summary, we have shown the following result.

\begin{theorem}\label{varinequalityV}
Let $V(t,S)$ be the price of an American style put option on underlying asset $S$ following a geometric L\'evy process with an admissible activity L\'evy measure $\nu$ satisfying the structural inequality (\ref{nuassumption}). Then $V$ is a solution to the linear complementarity problem:
\begin{eqnarray}
&&\partial_t V(t,S) + L^{S}[V](t,S)\le 0, \qquad V(t,S)\ge \Phi(S),
\label{varineq}
\\
&&
\left(\partial_t V(t,S) + L^{S}[V](t,S)\right)\cdot \left( V(t,S)- \Phi(S) \right) =0,
\label{complentarity}
\end{eqnarray}
for any $t\in[0,T), S>0$, and $V(T,S)=\Phi(S)=(K-S)^+$.
\end{theorem}

A standard method for solving the linear complementarity problem (\ref{varineq})--(\ref{complentarity}) is based on construction of an approximate solution by means of the penalty method. A nonnegative penalty function ${\mathcal G}_\varepsilon(t,V)$ penalizes negative values of the difference $ V(t,S) -\Phi(S)$. 
For example, one can consider the penalty function of the form:
\[
{\mathcal G}_\varepsilon(t,V)(S) = \varepsilon^{-1} \min(S/K,\, 1) ( \Phi(S) - V(t,S))^+,
\]
where $0<\varepsilon\ll 1$ is a small parameter. Clearly, ${\mathcal G}_\varepsilon(t,V)(S)=0$ if and only if $V(t,S)\ge \Phi(S)$. Then the penalized problem for the approximate solution $V=V_\varepsilon$ to (\ref{varineq})--(\ref{complentarity}) reads as follows:

\begin{eqnarray}
&&\partial_t V + L^{S}[V] + {\mathcal G}_\varepsilon(t,V)=0, 
\quad S>0, t\in[0,T),
\label{Vpenalized}
\\
&&
V(T,S) = \Phi(S).
\nonumber
\end{eqnarray}

In terms of the transformed function $u(\tau,x)=e^{r\tau} V(T-\tau, Ke^x)$ and the shifted function $U=u-u_{BS}$ the penalized PIDE problem (\ref{Vpenalized}) can be rewritten as follows:
\begin{equation}
\frac{\partial U}{\partial \tau} + A U = F(U) + h(\tau) + g_\varepsilon(\tau, U),
\quad U(0)=0.
\label{penalizedproblem_transformed-shifted}
\end{equation}

Equation (\ref{penalizedproblem_transformed-shifted}) can be understood as an abstract parabolic equation in the phase space $X^\gamma={\mathscr L}^p_{2\gamma}(\mathbb{R})$, i.e. $U(\tau)\in C([0,T], X^\gamma)$ where $F: X^\gamma\to X$. Furthermore,  $h(\tau), g_\varepsilon(\tau, U) \in X$ for any $\tau\in(0,T]$ and $U\in X^\gamma$, i.e. they are $x$-dependent functions for each $\tau$.

The penalty term $g_\varepsilon$ can be deduced from ${\mathcal G}_\varepsilon$, i.e. 
\[
g_\varepsilon(\tau,U(\tau,x))(x) = \varepsilon^{-1}   e^{x^-}
(w(\tau,x) - U(\tau,x))^+,
\quad \hbox{where}\ 
w(\tau,x) =   e^{r\tau}\Phi(K e^x) - u_{BS}(\tau,x).
\]
Recall that the linear operators $A$ and $f$ were defined in (\ref{Au_def}) and (\ref{functional_f_def}) and 
\[
F(U) = (r-\sigma^2/2)\frac{\partial U}{\partial x} + f[U],
\qquad h(\tau)=f[u_{BS}(\tau,\cdot)].
\]

Before proving existence and uniqueness of a solution to the penalized PIDE equation (\ref{penalizedproblem_transformed-shifted}) we need the following auxiliary lemma.

\begin{lemma}\label{gprop}
The penalty function $g_\varepsilon: [0,T]\times X \to X$ is Lipschitz continuous in the $U$ variable and H\"older continuous in the $\tau$ variable, i.e. there exists a constant $C_0>0$ such that 
\[
\Vert g_\varepsilon(\tau, U_1)-g_\varepsilon(\tau, U_2)\Vert_X
\le \varepsilon^{-1}\Vert U_1-U_2\Vert_X,\ \ 
\Vert g_\varepsilon(\tau_1, U)-g_\varepsilon(\tau_2, U)\Vert_X
\le \varepsilon^{-1} C_0 |\tau_1-\tau_2|^{\frac{p+1}{2p}}
\]
for any $U, U_1, U_2\in X$ and $\tau, \tau_1,\tau_2\in[0,T]$.
\end{lemma}

\noindent Proof. Note the inequality $|a^+-b^+|\le |a-b|$ for all $a,b\in\mathbb{R}$. As $e^{x^-}\le 1$, we obtain
\begin{eqnarray*}
\Vert g_\varepsilon(\tau, U_1)-g_\varepsilon(\tau, U_2)\Vert_{L^p}^p
&\le&\varepsilon^{-p}\int_{-\infty}^\infty \left| (w(\tau,x) - U_1(x))^+ - (w(\tau,x) - U_2(x))^+ \right|^p \ud x 
\\
&\le& \varepsilon^{-p} \int_{-\infty}^\infty | U_1(x)- U_2(x)|^p \ud x 
=\varepsilon^{-p} \Vert U_1-U_2\Vert_{L^p}^p.
\end{eqnarray*}
Moreover, it is easy to verify that the function $ e^{x^-} w(\tau,x)$ belongs to $X=L^p$ and 
\[
w(\tau,x) =  e^{r\tau}\Phi(K e^x) - K N(-d_2(\tau,x)) + K e^{r\tau +x} N(-d_1(\tau,x)).
\]
Hence $g_\varepsilon(\tau,0)\in X=L^p$ and $g_\varepsilon(\tau,\cdot):X \to X$ is well defined and Lipschitz continuous mapping for any $\tau\in[0,T]$.

Recall that $d_1-d_2 = \sigma \sqrt{\tau}$, 
$d_1+d_2 = 2(x+ r\tau)/\sigma \sqrt{\tau}$, and, consequently,
$e^{r\tau +x} N^\prime(-d_1) - N^\prime(-d_2) =0$. Since $N(-d_1)=1-N(d_1)$ we obtain
\begin{eqnarray*}
\partial_\tau w &=& 
  r e^{r\tau}\Phi(K e^x) +  r K e^{r\tau +x} N(-d_1) 
 - K N^\prime(-d_2) \frac{\sigma}{2\sqrt{\tau}} 
 \\
&=&  r e^{r\tau}\Phi(K e^x) - r v 
 - K \frac{e^{-d_2^2/2}}{\sqrt{2\pi}} \frac{\sigma}{2\sqrt{\tau}} 
\end{eqnarray*}
where the auxiliary function $v$ was defined as in (\ref{vtaux}). Therefore
\begin{eqnarray*}
\Vert e^{x^-} \partial_\tau w\Vert_{L^p} 
&\le& 
r e^{r\tau} \Vert e^{x^-}\Phi(K e^x)\Vert_{L^p} 
+r \Vert e^{x^-} v\Vert_{L^p} 
+ \frac{K \sigma}{2\sqrt{\tau}}
\left(
\int_{-\infty}^\infty e^{p x^-}
\frac{e^{-p d_2^2/2}}{(2\pi)^{p/2}} \ud x
\right)^{1/p}
\\
&\le& 
r K e^{r\tau} \Vert e^{x^-} 1_{x\le 0}\Vert_{L^p} 
+r \Vert v\Vert_{L^p} 
+ \frac{K \sigma}{2\sqrt{\tau}}
\left(
\int_{-\infty}^\infty
\frac{e^{-p \xi^2/2}}{(2\pi)^{p/2}} \sigma\sqrt{\tau} \ud \xi
\right)^{1/p}
\\
&\le& C_0 \tau^{\frac{1}{2p} -\frac12},
\end{eqnarray*}
where $C_0>0$ is a constant independent of $\tau\in(0,T]$. 
Thus
\begin{eqnarray*}
\Vert g_\varepsilon(\tau_1, U)-g_\varepsilon(\tau_2, U)\Vert_{L^p}^p
&=&\varepsilon^{-p}\int_{-\infty}^\infty  e^{p x^-}\left| (w(\tau_1,x) - U(x))^+ - (w(\tau_2,x) - U(x))^+ \right|^p \ud x 
\\
&\le& \varepsilon^{-p} \int_{-\infty}^\infty  
e^{p x^-}| w(\tau_1,x) - w(\tau_2,x) |^p \ud x 
\\
&=& \varepsilon^{-p} \Vert  e^{x^-}(w(\tau_1,\cdot )-w(\tau_2,\cdot ))\Vert_{L^p}^p.
\end{eqnarray*}
Hence 
\[
\Vert g_\varepsilon(\tau_1, U)-g_\varepsilon(\tau_2, U)\Vert_{L^p}
\le \varepsilon^{-1} \int_{\tau_1}^{\tau_2} \Vert e^{x^-} \partial_\tau w(\tau,\cdot)\Vert_{L^p} \ud\tau \le \varepsilon^{-1} C_0 |\tau_1-\tau_2|^{\frac{p+1}{2p}},
\]
as claimed. The proof of lemma follows.\hfill $\diamondsuit$

\medskip

Similarly as in the case of a linear PIDE, applying Lemmas~\ref{lemma-f}, \ref{lemmafuBS}, \ref{A-sectorial}, and Proposition~\ref{semilinear_existence_result} we obtain the following existence and uniqueness result for the nonlinear penalized PIDE (\ref{penalizedproblem_transformed-shifted}).

\begin{theorem}\label{existence_nonlinear_PIDE}
Assume $\nu$ is an admissible activity L\'evy  measure with the shape parameters $\alpha<3$, and either $\mu>0, D^\pm\in\mathbb{R}$, or $\mu=0$ and $D^- +1<0<D^+$. Let $X^\gamma={\mathscr L}^p_{2\gamma}(\mathbb{R})$ be the space of Bessel potentials where  $\frac12 \le \gamma<1$ and $\frac{\alpha-1}{2}<\gamma < \frac{p+1}{2p}$. Suppose that the structural condition (\ref{nuassumption}) is fulfilled for the L\'evy measure $\nu$. 

Then, for any $\varepsilon >0$ and $T>0$,  the nonlinear penalized PIDE (\ref{penalizedproblem_transformed-shifted}) has the unique solution $U_\varepsilon\in C([0,T),X^{\gamma})$. Moreover, $U_\varepsilon(\tau, \cdot )\in X^1={\mathscr L}^p_{2}(\mathbb{R})\hookrightarrow W^{2,p}(\mathbb{R})$, and $\partial_\tau U_\varepsilon(\tau, \cdot )\in L^p(\R)$ for any $\tau\in (0,T)$. 
\end{theorem}

\section{Numerical experiments}

In this section we present comparison of solutions to the linear PIDE with various L\'evy measures. We consider European style of put options only, i.e. $\Phi(S)=(K-S)^+$. We compare a solution for the linear Black-Scholes equation with solutions to the Merton and Variance Gamma PIDE models. The common model parameters were chosen as follows $\sigma=0.23, K=100, T=1$ and $r\in\{0, 0.1\}$. As for the underlying L\'evy process we consider the Variance Gamma process with parameters $\theta=-0.43, \kappa=0.27$ and the Merton processes with parameters $\lambda=0.1, m=-0.2, \delta=0.15$. In order to compute numerical solution we chose the finite difference discretization scheme proposed and analyzed by Cruz and \v{S}ev\v{c}ovi\v{c} in \cite{NBS19}. The scheme is based on a uniform spatial finite difference discretization with a spatial step $\Delta x = 0.01$, and implicit time discretization with a step $\Delta t = 0.005$. The total number of spatial discretization steps was chosen $N=400$ and the number of time discretization steps $M=200$. We restricted the spatial computational domain to $x\in[-L,L ]$ where $L=4$. We refer the reader to \cite{NBS19} for details concerning discretization scheme.

\begin{figure}[!ht]
\centering
    \includegraphics[width=0.48\textwidth]{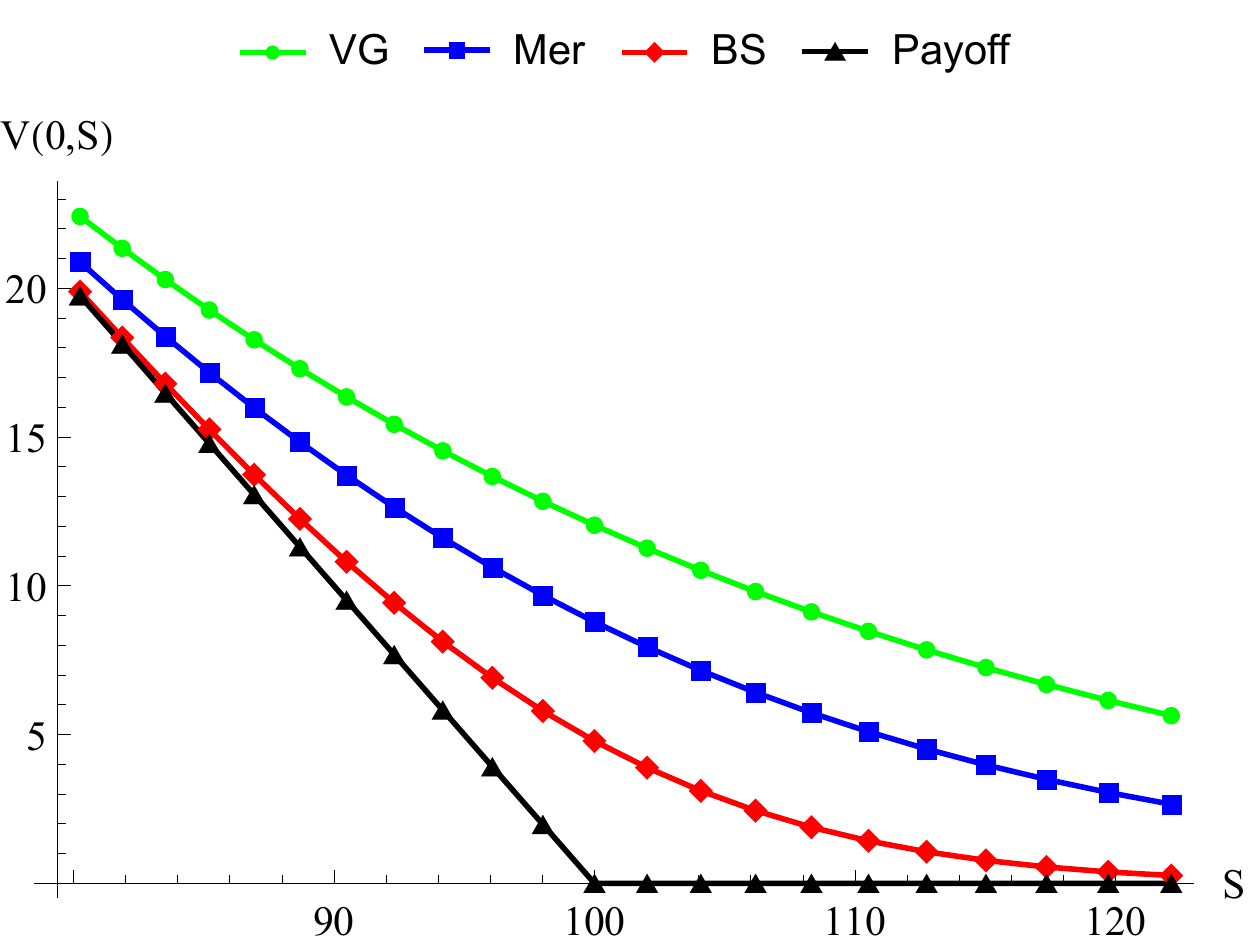}
\quad
    \includegraphics[width=0.48\textwidth]{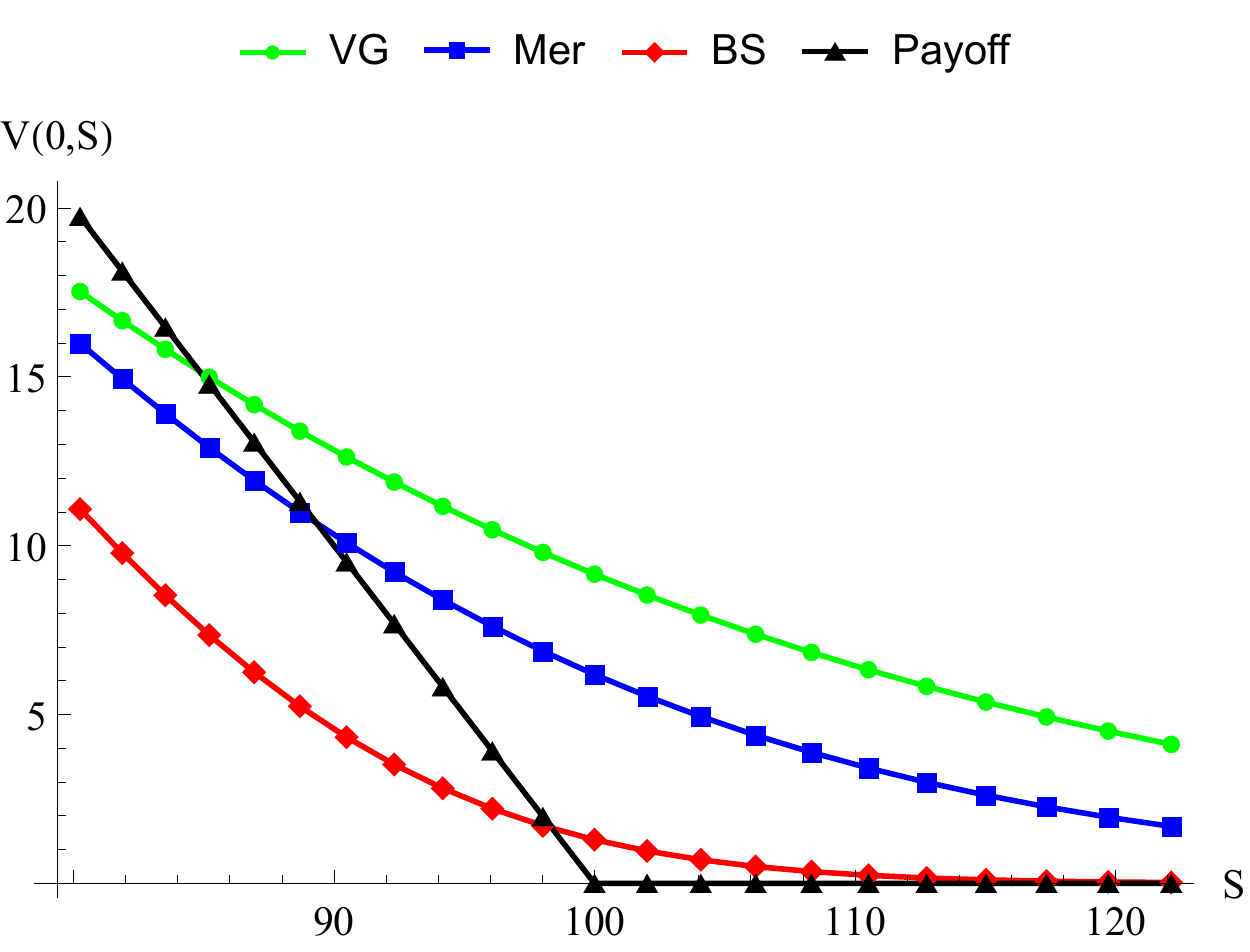}

    a) \hskip 7truecm b)

\caption{
Graphical comparison of European put option prices for the Black--Scholes (BS) model and the  PIDE Variance Gamma (VG) and Merton's (Mer) models.
} 
\label{fig-1}
\end{figure}

In Fig.~\ref{fig-1} we show comparison of European put option prices between PIDE models and the linear Black--Scholes model. In Fig.~\ref{fig-1} a) we plot put option prices $V(0,S)$ for $S\in[80,125]$ for the zero interest rate $r=0$, whereas b) depicts put option prices for the interest rate $r=0.1$. Numerical values of option prices are summarized in Table~\ref{tab1} for two different values of the interest rate $r=0.1$ and $r=0$. The option price for both Merton's as well as the Variance Gamma models are higher when compared to the option prices computed by means of the classical Black-Scholes model. This is in accordance with an intuitive observation that prices of put or call options should be higher on underlyings assets following stochastic processes with jumps when compared to those following a continuous geometric Brownian motion.

\begin{table}

\centering
\caption{European put option prices $V(0,S)$ for the Black-Scholes and PIDE models under Variance Gamma and Merton's processes for $r=0$ and $r=0.1$.}
\label{tab1}
\small

\medskip
\begin{tabular}{l|ll|ll|ll|l}
         & \multicolumn{2}{c|}{BS} &    \multicolumn{2}{c|}{PIDE-VG} & \multicolumn{2}{c|}{PIDE-Merton}  &  Payoff  \\ 
 $S$     &  $r=0$  & $r=0.1$ & $r=0$   & $r=0.1$ & $r=0$   & $r=0.1$ & \\
\hline\hline
 85.2144 & 15.2547 & 7.35166 & 19.2687 & 14.9855 & 17.1692 & 12.9056 & 14.7856\\
 88.692  & 12.2484 & 5.24145 & 17.2948 & 13.3899 & 14.8335 & 10.9901 & 11.308\\
 92.3116 & 9.42895 & 3.51944 & 15.428  & 11.8822 & 12.6423 & 9.21922 & 7.68837\\
 96.0789 & 6.90902 & 2.21106 & 13.674  & 10.4691 & 10.6201 & 7.61307 & 3.92106\\
 100.    & 4.78444 & 1.29196 & 12.0372 & 9.15576 & 8.78655 & 6.18483 & 0.\\
 104.081 & 3.1099  & 0.69843 & 10.52   & 7.94499 & 7.155   & 4.94044 & 0.\\
 108.329 & 1.88555 & 0.34773 & 9.12343 & 6.83762 & 5.73137 & 3.87864 & 0.\\
 112.75  & 1.0604  & 0.15881 & 7.84623 & 4.51403 & 5.83246 & 2.99166 & 0.\\
\hline
\end{tabular}

\end{table}

\section{Conclusions}
In this paper, we analyzed existence and uniqueness of solutions to a partial integro-differential equation (PIDE) in the Bessel potential space. As a motivation we considered a model for pricing vanilla call and put options on underlying assets following a geometric L\'evy stochastic process. Using the theory of abstract semilinear parabolic equations we proved existence and uniqueness of solutions in the Bessel potential space representing a fractional power space of the space of Lebesgue $p$-integrable functions with respect to the second order Laplace differential operator. We generalized known existence results for a wider class of L\'evy measures including those having strong singular kernel. We also proved existence and uniqueness of solutions to the penalized PIDE representing approximation of the linear complementarity problem arising in pricing American style of options. 

\section{Acknowledgements}

This research was supported by the project CEMAPRE MULTI/00491 financed by FCT/MEC through national funds and the Slovak research Agency Project VEGA 1/0062/18.

\bibliographystyle{apa}

\bibliography{paper}

\end{document}